\documentclass[a4paper,11pt]{amsart}
\usepackage{graphicx}
 \usepackage{mathptmx}
\usepackage{amsmath}
\usepackage{amssymb}
\usepackage{enumitem}
\usepackage{xcolor}

\newmuskip\pFqmuskip

\newcommand*\pFq[6][8]{%
  \begingroup 
  \pFqmuskip=#1mu\relax
  \mathcode`=\string"8000
  \begingroup\lccode`\~=`\,
  \lowercase{\endgroup\let~}\pFqcomma
  F^{#2}_{#3}{\left(\genfrac..{0pt}{}{#4}{#5}\bigg|#6\right)}%
  \endgroup
}
\newcommand{\pFqcomma}{\mskip\pFqmuskip}

\newtheorem{theorem}{Theorem}

\newtheorem{corollary}[theorem]{Corollary}

\begin{document}

\title[]{Some identities on degenerate harmonic and degenerate higher-order harmonic numbers}

\author{Taekyun  Kim}
\address{Department of Mathematics, Kwangwoon University, Seoul 139-701, Republic of Korea}
\email{tkkim@kw.ac.kr}

\author{Dae San  Kim }
\address{Department of Mathematics, Sogang University, Seoul 121-742, Republic of Korea}
\email{dskim@sogang.ac.kr}


\begin{abstract}
The harmonic numbers and higher-order harmonic numbers appear frequently in several areas which are related to combinatorial identities, many expressions involving special functions in analytic number theory, and analysis of algorithms. The aim of this paper is to study the degenerate harmonic and degenerate higher-order harmonic numbers, which are respectively degenerate versions of the harmonic and higher-order harmonic numbers, in connection with the degenerate zeta and degenerate Hurwitz zeta function. Here the degenerate zeta and degenerate Hurwitz zeta function are respectively degenerate versions of the Riemann zeta and Hurwitz zeta function. We show that several infinite sums involving the degenerate higher-order harmonic numbers can be expressed in terms of the degenerate zeta function. Furthermore, we demonstrate that an infinite sum involving finite sums of products of the degenerate harmonic numbers can be represented by using the degenerate Hurwitz zeta function.
\end{abstract}

\maketitle

\markboth{\centerline{\scriptsize Some identities on degenerate harmonic and degenerate higher-order harmonic numbers}}
{\centerline{\scriptsize T. Kim and D. S. Kim}}
 
\section{Introduction} 
We have witnessed that explorations for degenerate versions of some special numbers and polynomials regained interests of some mathematicians and produced many fascinating results, which began with the pioneering work of Carlitz on the degenerate Bernoulli and degenerate Euler polynomials. These have been done by using many tools like combinatorial methods, generating functions, umbral calculus, $p$-adic analysis, probability theory, special functions, quantum mechanics, analytic number theory, operator theory and differential equations. \par
The aim of this paper is to study the degenerate harmonic and degenerate higher-order harmonic numbers in connection with the degenerate zeta and degenerate Hurwitz zeta function. Here the degenerate harmonic numbers, the degenerate higher-order harmonic numbers, the degenerate zeta function and the degenerate Hurwitz zeta function are respectively degenerate versions of the harmonic numbers, the higher-order harmonic numbers, Riemann zeta function and Hurwitz zeta function. We show that several infinite sums involving the degenerate higher-order harmonic numbers can be expressed in terms of the degenerate zeta function. Furthermore, we demonstrate that an infinite sum involving finite sums of products of the degenerate harmonic numbers can be represented by using the degenerate Hurwitz zeta function. \par
The outline of this paper is as follows. In Section 1, we recall the harmonic numbers and the generalized harmonic nmubers of order $\alpha$. We remind the reader of the degenerate exponentials, the degenerate logarithms and of the degenerate polylogarithms $\mathrm{Li}_{k,\lambda}(t)$. We recall the degenerate harmonic numbers, the degenerate higher-order harmonic numbers $H_{k,\lambda}^{(p)}$ and the degenerate zeta function $\zeta_{\lambda}(s)$. Section 2 is the main result of this paper. In Theorem 1, we evaluate the integral $I_{\lambda}(r,p)=\int_{0}^{1}x^{r-1}\mathrm{Li}_{p,\lambda}(x)dx$, for any positive integers $r,p$ with $p \ge 2$, as the sum of two finite sums, one of which involves the degenerate zeta function. In Theorem 2, we express the infinite sum $\sum_{k=1}^{\infty}\frac{H_{k,\lambda}^{(p)}}{\binom{k+n-1}{n}(k+n)}$, as a finite sum involving $I_{\lambda}(j+1,p+1)$, for $j=0,1,\dots,n-2$. In theorem 3, we represent the infinite sum $\sum_{k=1}^{\infty}\frac{(-1)^{k-1}\lambda^{k-1}(1)_{k,1/\lambda}}{(k-1)!k^{p}(k+1)}$, as a finite sum involving the degenerate zeta function. In Theorem 4, we express the infinite sum $\sum_{k=1}^{\infty}\frac{nH_{k,\lambda}^{(p)}}{\binom{k+n-1}{n}(k+n)}$, as another infinite sum not involving the degenerate higher-order harmonic numbers. In Theorem 5, we show that $\sum_{k=1}^{\infty}\frac{H_{k,\lambda}^{(p)}}{k(k+1)}=\zeta_{\lambda}(p+1)$, and that $\sum_{k=1}^{\infty}\frac{H_{k,\lambda}^{(p)}}{k(k+2)}$ is equal to a finite sum involving the degenerate zeta function. Finally, in Corollary 7 we show that $\sum_{k=p}^{\infty}\frac{H_{k,p,\lambda}}{\binom{k+n-1}{n}(n+k)}=p!\zeta_{\lambda}^{*}(p+1,n-p+1)$, where $\zeta_{\lambda}^{*}(k,x)$ is the degenerate Hurwitz zeta function and $H_{k,p,\lambda}=\sum_{k_{1}+\cdots+k_{p}=k}H_{k_{1},\lambda}H_{k_{2},\lambda}\cdots H_{k_{p},\lambda},\,\,(k_{1},\dots,k_{p} \ge 1)$. For the rest of this section, we recall the facts that are needed throughout this paper.

\vspace{0.1in}

It is well known that the harmonic numbers are defined by 
\begin{equation}
-\frac{\log (1-t)}{1-t}=\sum_{n=1}^{\infty}H_{n}t^{n},\quad (\mathrm{see}\ [6,7,11,12,22]). \label{1}
\end{equation}
By \eqref{1}, we get $H_{0}=1,\ H_{n}=1+\frac{1}{2}+\cdots+\frac{1}{n},\ (n\ge 1)$. \\
For $\alpha\in\mathbb{N}$, the generalized harmonic numbers of order $\alpha$ are defined by 
\begin{equation}
H_{0}^{(\alpha)}=0,\quad H_{n}^{(\alpha)}=1+\frac{1}{2^{\alpha}}+\frac{1}{3^{\alpha}}+\cdots+\frac{1}{n^{\alpha}},\quad (n\ge 1),\quad (\mathrm{see}\ [12,21,22]). \label{2}	
\end{equation} \par
For $s\in\mathbb{C}$ with $\mathrm{Re}(s)>0$, the gamma function is defined by 
\begin{equation}
\Gamma(s)=\int_{0}^{\infty}e^{-t}t^{s-1}dt,\quad (\mathrm{see}\ [1-25]). \label{4}
\end{equation}
From \eqref{4}, we note that $\Gamma(s+1)=s\Gamma(s)$. For $\mathrm{Re}(\alpha)>0$ and $\mathrm{Re}(\beta)>0$, the beta function is defined by 
\begin{equation}
B(\alpha,\beta)=\int_{0}^{1}t^{\alpha-1}(1-t)^{\beta-1}dt,\quad (\mathrm{see}\ [1-25]). \label{5}
\end{equation}
By \eqref{4} and \eqref{5}, we get $B(\alpha,\beta)=\frac{\Gamma(\alpha)\Gamma(\beta)}{\Gamma(\alpha+\beta)}$. \par
Let $\lambda$ be any nonzero real number. Then the degenerate exponentials are given by 
\begin{equation}
e_{\lambda}^{x}(t)=\sum_{k=0}^{\infty}\frac{(x)_{k,\lambda}}{k!}t^{k},\quad (\mathrm{see}\ [7,10]),\label{6}
\end{equation}
where 
\begin{equation}
(x)_{0,\lambda}=1,\quad (x)_{n,\lambda}=x(x-\lambda)(x-2\lambda)\cdots(x-(n-1)\lambda),\quad (n\ge 1). \label{7}
\end{equation}
When $x=1$, $e_{\lambda}(t)=e_{\lambda}^{1}(t)=\sum_{k=0}^{\infty}\frac{(1)_{k,\lambda}}{k!}t^{k}=(1+\lambda t)^{\frac{1}{\lambda}},$ (see [8,11]). Let $\log_{\lambda}(t)$ be the degenerate logarithm, which is the compositional inverse of $e_{\lambda}(t)$, so that $\log_{\lambda}(e_{\lambda}(t))=e_{\lambda}(\log_{\lambda}(t))=t$. \\
Then we have 
\begin{equation}  
\log_{\lambda}(1+t)=\sum_{k=1}^{\infty}\frac{\lambda^{k-1}(1)_{k,1/\lambda}}{k!}t^{k},\quad (|t|<1),\quad (\mathrm{see}\ [7,14]). \label{8}
\end{equation}
Note that $\lim_{\lambda\rightarrow 0}\log_{\lambda}(1+t)=\log(1+t)$. \\
The degenerate polylogarithms are defined by Kim-Kim as
\begin{equation}
\mathrm{Li}_{k,\lambda}(t)=\sum_{n=1}^{\infty}\frac{(-1)^{n-1}\lambda^{n-1}(1)_{n,1/\lambda}}{(n-1)!n^{k}}t^{n},\quad (k\in\mathbb{Z},\ |t|<1),\quad (\mathrm{see}\ [7]). \label{9}	
\end{equation}
Note that $$\mathrm{Li}_{1,\lambda}(t)=-\log_{\lambda}(1-t),\quad \lim_{\lambda\rightarrow 0}\mathrm{Li}_{k,\lambda}(t)=\mathrm{Li}_{k}(t)=\sum_{n=1}^{\infty}\frac{t^{n}}{n^{k}},$$
where $\mathrm{Li}_{k}(t)$ is the polylogarithm of order $k$. \par
In [12], the degenerate harmonic numbers are given by 
\begin{equation}
-\frac{\log_{\lambda}(1-t)}{1-t}=\sum_{n=0}^{\infty}H_{n,\lambda}t^{n},\quad (|t|<1). \label{10}
\end{equation}
Thus, by \eqref{10}, we get 
\begin{equation}
H_{0,\lambda}=0,\quad H_{n,\lambda}=\sum_{k=1}^{n}\frac{\lambda^{k-1}(-1)^{k-1}(1)_{k,1/\lambda}}{k!},\quad(n \ge 1).\label{11}	
\end{equation}
Note that 
$$\lim_{\lambda\rightarrow 0}H_{n,\lambda}=H_{n},\quad (n\ge 1).$$ \\
For $s\in\mathbb{C}$ with $\mathrm{Re}(s)>1$, the degenerate zeta function is defined by 
\begin{equation} 
\zeta_{\lambda}(s)=\sum_{n=1}^{\infty}\frac{(-1)^{n-1}\lambda^{n-1}(1)_{n,1/\lambda}}{(n-1)!n^{s}},\quad (\mathrm{see}\ [12]). \label{12}
\end{equation}
Note that $\lim_{\lambda\rightarrow 0}\zeta_{\lambda}(s)=\sum_{n=1}^{\infty}\frac{1}{n^{s}}=\zeta(s)$ is the Riemann zeta function.
The degenerate higher-order harmonic numbers are defined by 
\begin{equation}
\frac{\mathrm{Li}_{k,\lambda}(t)}{1-t}=\sum_{n=1}^{\infty}H_{n,\lambda}^{(k)}t^{n},\quad (|t|<1),\quad (\mathrm{see}\ [12]). \label{13}	
\end{equation}
Thus, by \eqref{13}, we get 
\begin{equation}
H_{n,\lambda}^{(k)}=\sum_{l=1}^{n}\frac{(-\lambda)^{l-1}(1)_{l,1/\lambda}}{l^{k}(l-1)!},\quad (n\ge 1),\quad H_{0,\lambda}^{(k)}=0.\label{14}
\end{equation}
Note that 
\begin{displaymath}
	\lim_{\lambda\rightarrow 0}H_{n,\lambda}^{(k)}=H_{n}^{(k)}=1+\frac{1}{2^{k}}+\frac{1}{3^{k}}+\cdots+\frac{1}{n^{k}}. 
\end{displaymath}

\section{Some identities on degenerate harmonic and degenerate higher-order harmonic numbers} 
For $r,p\in\mathbb{N}$ with $p>1$, we let 
\begin{equation}
I_{\lambda}(r,p)=\int_{0}^{1}x^{r-1}\mathrm{Li}_{p,\lambda}(x)dx. \label{16} 
\end{equation}
Then, by \eqref{9} and \eqref{16}, we get 
\begin{align}
I_{\lambda}(r,p)&=\int_{0}^{1}x^{r-1}\mathrm{Li}_{p,\lambda}(x)dx=\sum_{n=1}^{\infty}\frac{(-\lambda)^{n-1}(1)_{n,1/\lambda}}{(n-1)!n^{p}}\int_{0}^{1}x^{n+r-1}dp \label{17} \\
&=\sum_{n=1}^{\infty}\frac{(-\lambda)^{n-1}(1)_{n,1/\lambda}}{(n-1)!n^{p}}\frac{1}{n+r}. \nonumber	
\end{align}
From \eqref{17}, we have 
\begin{align}
I_{\lambda}(r,p)&=\frac{1}{r}\sum_{n=1}^{\infty}\frac{(-\lambda)^{n-1}(1)_{n,1/\lambda}}{(n-1)!n^{p}}\frac{n+r-n}{n+r} \label{18} \\
&=\frac{1}{r}\sum_{n=1}^{\infty}\frac{(-\lambda)^{n-1}(1)_{n,1/\lambda}}{(n-1)!n^{p}}-\frac{1}{r}\sum_{n=1}^{\infty}\frac{(-\lambda)^{n-1}(1)_{n,1/\lambda}}{(n-1)!n^{p-1}}\frac{1}{n+r} \nonumber \\
&=\frac{1}{r}\zeta_{\lambda}(p)-\frac{1}{r}I_{\lambda}(r,p-1). \nonumber	
\end{align}
By \eqref{18}, we get 
\begin{align}
I_{\lambda}(r,p)&=\frac{1}{r}\zeta_{\lambda}(p)-\frac{1}{r}I_{\lambda}(r,p-1) \label{19} \\
&=\frac{1}{r}\zeta_{\lambda}(p)-\frac{1}{r^{2}}\zeta_{\lambda}(p-1)+\frac{1}{r^{2}}I_{\lambda}(r,p-2)\nonumber\\
&=\cdots\nonumber \\
&=\sum_{k=1}^{p-1}\frac{(-1)^{k-1}}{r^{k}}\zeta_{\lambda}(p-k+1)+\frac{(-1)^{p-1}}{r^{p-1}}I_{\lambda}(r,1). \nonumber	
\end{align}
Now, we observe that 
\begin{align}
I_{\lambda}(r,1)&=\int_{0}^{1}x^{r-1}\mathrm{Li}_{1,\lambda}(x)dx=-\int_{0}^{1}x^{r-1}\log_{\lambda}(1-x)dx \label{20}\\
&=\bigg[\frac{1-x^{r}}{r}\log_{\lambda}(1-x)\bigg]_{0}^{1}+\frac{1}{r}\int_{0}^{1}(1-x^{r})(1-x)^{\lambda-1}dx\nonumber \\
&=\frac{1}{r}\sum_{i=0}^{r-1}\int_{0}^{1}x^{i}(1-x)^{\lambda}dx. \nonumber
\end{align}
Integrating by parts, we get 
\begin{align}
&\int_{0}^{1}x^{i}(1-x)^{\lambda}dx=\bigg[-\frac{1}{\lambda+1}(1-x)^{\lambda+1}x^{i}\bigg]_{0}^{1}+\frac{i}{\lambda +1}\int_{0}^{1}x^{i-1}(1-x)^{\lambda+1}dx \label{21} \\
&=\frac{i}{\lambda+1}\int_{0}^{1}x^{i-1}(1-x)^{\lambda+1}dx=\frac{i(i-1)}{(\lambda+1)(\lambda+2)}\int_{0}^{1}x^{i-2}(1-x)^{\lambda+2}dx=\cdots \nonumber \\
&=\frac{i(i-1)\cdots 1}{(\lambda+1)(\lambda+2)\cdots(\lambda+i)}\int_{0}^{1}(1-x)^{\lambda+i}dx=\frac{i!}{(\lambda+1)(\lambda+2)\cdots(\lambda+i)}\frac{1}{\lambda+i+1} \nonumber \\
&=\frac{1}{\binom{\lambda+i}{i}(\lambda+i+1)}.\nonumber
\end{align}
From \eqref{20} and \eqref{21}, we have 
\begin{equation}
I_{\lambda}(r,1)=\frac{1}{r}\sum_{i=0}^{r-1}\frac{1}{\binom{\lambda+i}{i}(\lambda+i+1)}. \label{22}
\end{equation}
By \eqref{19} and \eqref{22}, we get 
\begin{equation} 
I_{\lambda}(r,p)=\sum_{k=1}^{p-1}\frac{(-1)^{k-1}}{r^{k}}\zeta_{\lambda}(p-k+1)+\frac{(-1)^{p-1}}{r^{p}}\sum_{i=0}^{r-1}\frac{1}{\binom{\lambda+i}{i}(\lambda+i+1)}. \label{23}	
\end{equation}
Therefore, by \eqref{23}, we obtain the following theorem. 
\begin{theorem}
For $p,r\in\mathbb{N}$ with $p>1$, we have 
\begin{align*}
I_{\lambda}(r,p)&=\int_{0}^{1}x^{r-1}\mathrm{Li}_{p,\lambda}(x)dx \\
&=\sum_{k=1}^{p-1}\frac{(-1)^{k-1}}{r^{k}}\zeta_{\lambda}(p-k+1)+\frac{(-1)^{p-1}}{r^{p}}\sum_{i=0}^{r-1}\frac{1}{\binom{\lambda+i}{i}(\lambda+i+1)}. 
\end{align*}
\end{theorem}
From \eqref{9}, we note that 
\begin{equation}
\frac{d}{dx}\mathrm{Li}_{p,\lambda}(x)=\frac{1}{x}\sum_{n=1}^{\infty}\frac{(-\lambda)^{n-1}(1)_{n,1/\lambda}}{(n-1)!n^{p-1}}x^{n}=\frac{1}{x}\mathrm{Li}_{p-1,\lambda}(x). \label{24}
\end{equation}
By \eqref{5}, we get 
\begin{align}
&n!\sum_{k=1}^{\infty}\frac{H_{k,\lambda}^{(p)}}{k(k+1)\cdots(k+n)}=\sum_{k=1}^{\infty}H_{k,\lambda}^{(p)}\frac{\Gamma(k)\Gamma(n+1)}{\Gamma(k+n+1)}
\label{25} \\
&=\sum_{k=1}^{\infty}H_{k,\lambda}^{(p)}B(k,n+1)=\sum_{k=1}^{\infty}H_{k,\lambda}^{(p)}\int_{0}^{1}t^{k-1}(1-t)^{n}dt \nonumber\\
&=\int_{0}^{1}\bigg(\sum_{k=1}^{\infty}H_{k,\lambda}^{(p)}t^{k}\bigg)\frac{(1-t)^{n}}{t}dt =\int_{0}^{1}\frac{\mathrm{Li}_{p,\lambda}(t)}{t}(1-t)^{n-1}dt\nonumber \\
&=\int_{0}^{1}\bigg(\frac{d}{dt}\mathrm{Li}_{p+1,\lambda}(t)\bigg)(1-t)^{n-1}dt=\Big[(1-t)^{n-1}\mathrm{Li}_{p+1,\lambda}(t)\Big]_{0}^{1}+(n-1)\int_{0}^{1}\mathrm{Li}_{p+1,\lambda}(t)(1-t)^{n-2}dt \nonumber \\
&=(n-1)\sum_{j=0}^{n-2}\binom{n-2}{j}(-1)^{j}\int_{0}^{1}t^{j}\mathrm{Li}_{p+1,\lambda}(t)dt=(n-1)\sum_{j=0}^{n-2}\binom{n-2}{j}(-1)^{j}I_{\lambda}(j+1,p+1). \nonumber
\end{align}
Thus, by \eqref{25}, we get 
\begin{equation}
\sum_{k=1}^{\infty}\frac{n!H_{k,\lambda}^{(p)}}{k(k+1)\cdots(k+n)}=(n-1)\sum_{j=0}^{n-2}\binom{n-2}{j}(-1)^{j}I_{\lambda}(j+1,p+1).\label{26}
\end{equation}
In particular, from Theorem 1, for $n=2$ we have 
\begin{equation}
2!\sum_{k=1}^{\infty}\frac{H_{k,\lambda}^{(p)}}{k(k+1)(k+2)}=I_{\lambda}(1,p+1)=\sum_{k=1}^{p}(-1)^{k-1}\zeta_{\lambda}(p-k+2)+\frac{(-1)^{p}}{\lambda+1}. \label{27}	
\end{equation}
Therefore, by \eqref{26} and \eqref{27}, we obtain the following theorem. 
\begin{theorem}
For any integer $p\ge 1$, we have 
\begin{displaymath}
\sum_{k=1}^{\infty}\frac{H_{k,\lambda}^{(p)}}{k(k+1)(k+2)}=\frac{1}{2}\bigg(\sum_{k=1}^{p}(-1)^{k-1}\zeta_{\lambda}(p+2-k)+\frac{(-1)^{p}}{\lambda+1}\bigg). 
\end{displaymath}
More generally, for any integers $n,p$, with $n \ge 2, p \ge1$, we have 
\begin{displaymath}
\sum_{k=1}^{\infty}\frac{H_{k,\lambda}^{(p)}}{\binom{k+n-1}{n}(k+n)}=(n-1)\sum_{j=0}^{n-2}\binom{n-2}{j}(-1)^{j}I_{\lambda}(j+1,p+1). 
\end{displaymath}
\end{theorem}
From \eqref{9}, we note that 
\begin{align}
&(n-1)!\sum_{k=1}^{\infty}\frac{(-\lambda)^{k-1}(1)_{k,1/\lambda}}{(k-1)!k^{p}(k+1)\cdots (k+n)}==\sum_{k=1}^{\infty}\frac{(-\lambda)^{k-1}(1)_{k,1/\lambda}}{k^{p}(p-1)!}B(k+1,n) \label{28} \\
&=\sum_{k=1}^{\infty}\frac{(-\lambda)^{k-1}(1)_{k,1/\lambda}}{k^{p}(k-1)!}\int_{0}^{1}t^{k}(1-t)^{n-1}dt=\int_{0}^{1}\bigg(\sum_{k=1}^{\infty}\frac{(-\lambda)^{k-1}(1)_{k,1/\lambda}}{k^{p}(k-1)!}t^{k}\bigg)(1-t)^{n-1}dt \nonumber \\
&=\int_{0}^{1}\mathrm{Li}_{p,\lambda}(t)(1-t)^{n-1}dt=\sum_{j=0}^{n-1}\binom{n-1}{j}(-1)^{j}\int_{0}^{1}\mathrm{Li}_{p,\lambda}t^{j}dt \nonumber\\
&=\sum_{j=0}^{n-1}\binom{n-1}{j}(-1)^{j}I_{\lambda}(j+1,p).\nonumber	
\end{align}
Thus, by \eqref{28}, we get 
\begin{equation}
(n-1)!\sum_{k=1}^{\infty}\frac{(-\lambda)^{k-1}(1)_{k,1/\lambda}}{(k-1)!k^{p}(k+1)\cdots(k+n)}=\sum_{j=0}^{n-1}\binom{n-1}{j}(-1)^{j}I_{\lambda}(j+1,p), \label{29}
\end{equation}
where $n$ is a positive integer, \par 
Taking $n=1$ in \eqref{29}, we have 
\begin{align}
\sum_{k=1}^{\infty}\frac{(-1)^{k-1}\lambda^{k-1}(1)_{k,1/\lambda}}{(k-1)!k^{p}(k+1)}&=I_{\lambda}(1,p)=\sum_{k=1}^{p-1}(-1)^{k-1}\zeta_{\lambda}(p-k+1)+\frac{(-1)^{p-1}}{\lambda+1}. 	\label{30}
\end{align}
Therefore, by \eqref{30}, we obtain the following theorem. 
\begin{theorem}
For any integer $p$, with $p\ge 2$, we have 
\begin{displaymath}
\sum_{k=1}^{\infty}\frac{(-1)^{k-1}\lambda^{k-1}(1)_{k,1/\lambda}}{(k-1)!k^{p}(k+1)}=\sum_{k=1}^{p-1}(-1)^{k-1}\zeta_{\lambda}(p-k+1)+\frac{(-1)^{p-1}}{\lambda+1}. 
\end{displaymath}
\end{theorem}
Replacing $n$ by $n-1$ and $p$ by $p+1$ in \eqref{29}, we get 
\begin{equation}
(n-2)!\sum_{k=1}^{\infty}\frac{(-\lambda)^{k-1}(1)_{k,1/\lambda}}{(k-1)!k^{p+1}(k+1)\cdots(k+n-1)}=\sum_{j=0}^{n-2}\binom{n-2}{j}(-1)^{j}I_{\lambda}(j+1,p+1). \label{31}	
\end{equation}
Thus, by \eqref{26} and \eqref{31}, we get 
\begin{align}
&(n-1)!\sum_{k=1}^{\infty}\frac{(-\lambda)^{k-1}(1)_{k,1/\lambda}}{(k-1)!k^{p+1}(k+1)\cdots(k+n-1)}=(n-1)\sum_{j=0}^{n-2}\binom{n-2}{j}(-1)^{j}I_{\lambda}(j+1,p+1) \label{32} \\
&=\sum_{k=1}^{\infty}\frac{n!H_{k,\lambda}^{(p)}}{k(k+1)\cdots(k+n)}=\sum_{k=1}^{\infty}\frac{H_{k,\lambda}^{(p)}}{\binom{k+n-1}{n}(k+n)}. \nonumber
\end{align}
From \eqref{32}, we get 
\begin{align}
n\sum_{k=1}^{\infty}\frac{H_{k,\lambda}^{(p)}}{\binom{k+n-1}{n}(k+n)}&=n!\sum_{k=1}^{\infty}\frac{(-\lambda)^{k-1}(1)_{k,1/\lambda}}{(k-1)!k^{p+1}(k+1)\cdots(k+n-1)} \label{33} \\
&=\sum_{k=1}^{\infty}\frac{(-\lambda)^{k-1}(1)_{k,1/\lambda}}{(k-1)!k^{p}\binom{k+n-1}{n}}.\nonumber
\end{align}
Therefore, by \eqref{33}, we obtain the following theorem. 
\begin{theorem}
For any integer $n\ge 1$, we have 
\begin{displaymath}
\sum_{k=1}^{\infty}\frac{nH_{k,\lambda}^{(p)}}{\binom{k+n-1}{n}(k+n)}=\sum_{k=1}^{\infty}\frac{(-\lambda)^{k-1}(1)_{k,1/\lambda}}{(k-1)!k^{p}\binom{k+n-1}{n}}. 
\end{displaymath}
\end{theorem}
Let $m,n$ be integers with $m>n \ge 0$. Then, by \eqref{13}, we get 
\begin{align}
\int_{0}^{1}\frac{\mathrm{Li}_{p,\lambda}(x)}{x(1-x)}\big(1-x^{m-n}\big)x^{n}dx&=\sum_{k=1}^{\infty}H_{k,\lambda}^{(p)}\int_{0}^{1}x^{n+k-1}(1-x^{m-n})dx \label{34} \\
&=\sum_{k=1}^{\infty}H_{k,\lambda}^{(p)}\bigg(\frac{1}{n+k}-\frac{1}{m+k}\bigg) \nonumber \\
&=(m-n)\sum_{k=1}^{\infty}\frac{H_{k,\lambda}^{(p)}}{(n+k)(m+k)}. \nonumber	
\end{align}
From \eqref{34}, we note that 
\begin{align}
&\sum_{k=1}^{\infty}\frac{H_{k,\lambda}^{(p)}}{(n+k)(m+k)}=\frac{1}{m-n}\int_{0}^{1}\frac{\mathrm{Li}_{p,\lambda}(x)}{x(1-x)}(1-x^{m-n})x^{n}dx \label{35} \\
&=\frac{1}{m-n}\int_{0}^{1}\bigg(\frac{d}{dx}\mathrm{Li}_{p+1,\lambda}(x)\bigg)\frac{1-x^{m-n}}{1-x}x^{n}dx \nonumber  \\ 
&=\frac{1}{m-n}\int_{0}^{1}\bigg(\frac{d}{dx}\mathrm{Li}_{p+1,\lambda}(x)\bigg)(x^{n}+x^{n+1}+\cdots+x^{m-1})dx \nonumber \\
&=\frac{1}{m-n}\Big[\mathrm{Li}_{p+1,\lambda}(x)(x^{n}+x^{n+1}+\cdots+x^{m-1})\Big]_{0}^{1}-\frac{1}{m-n}\int_{0}^{1}\mathrm{Li}_{p+1,\lambda}(x)\frac{d}{dx}(x^{n}+\cdots+x^{m-1})dx \nonumber \\
&=\zeta_{\lambda}(p+1)-\frac{1}{m-n}\int_{0}^{1}\mathrm{Li
}_{p+1,\lambda}(x)\Big(nx^{n-1}+(n+1)x^{n}+\cdots+(m-1)x^{m-2}\Big)dx \nonumber  \\
&=\zeta_{\lambda}(p+1)-\frac{1}{m-n}\sum_{j=n}^{m-1}j\int_{0}^{1}\mathrm{Li}_{p+1,\lambda}(x)x^{j-1}dx=\zeta_{\lambda}(p+1)-\frac{1}{m-n}\sum_{j=n}^{m-1}jI_{\lambda}(j,p+1). \nonumber
\end{align}
Thus, by \eqref{35}, we get 
\begin{equation}
\sum_{k=1}^{\infty}\frac{H_{k,\lambda}^{(p)}}{(n+k)(m+k)}=\zeta_{\lambda}(p+1)-\frac{1}{m-n}\sum_{j=n}^{m-1}jI_{\lambda}(j,p+1). 	\label{36}
\end{equation}
In particular, for $m=1,\ n=0$, we have 
\begin{equation}
\sum_{k=1}^{\infty}\frac{H_{k,\lambda}^{(p)}}{k(k+1)}=\zeta_{\lambda}(p+1). \label{37}
\end{equation}
Taking $m=2$ and $n=0$ in \eqref{36}, from Theorem 1 we have 
\begin{align}
\sum_{k=1}^{\infty}\frac{H_{k,\lambda}^{(p)}}{k(k+2)}&=\zeta_{\lambda}(p+1)-\frac{1}{2}I_{\lambda}(1,p+1)\label{38} \\
&=\zeta_{\lambda}(p+1)-\frac{1}{2}\bigg(\sum_{k=1}^{p}(-1)^{k-1}\zeta_{\lambda}(p+2-k)+\frac{(-1)^{p}}{\lambda+1}\bigg) \nonumber \\
&=\frac{1}{2}\bigg(\zeta_{\lambda}(p+1)+ \zeta_{\lambda}(p)-\zeta_{\lambda}(p-1)+\cdots+(-1)^{p} \zeta_{\lambda}(2)+\frac{(-1)^{p+1}}{\lambda+1}\bigg).\nonumber
\end{align}
Therefore, by \eqref{37} and \eqref{38}, we obtain the following theorem. 
\begin{theorem}
For any integer $p\ge 1$, we have 
\begin{displaymath}
\sum_{k=1}^{\infty}\frac{H_{k,\lambda}^{(p)}}{k(k+1)}=\zeta_{\lambda}(p+1),
\end{displaymath}
and 
\begin{displaymath}
\sum_{k=1}^{\infty}\frac{H_{k,\lambda}^{(p)}}{k(k+2)}=\frac{1}{2}\Big(\zeta_{\lambda}(p+1)+ \zeta_{\lambda}(p)-\zeta_{\lambda}(p-1)+\cdots+(-1)^{p} \zeta_{\lambda}(2)+\frac{(-1)^{p+1}}{\lambda+1}\bigg).
\end{displaymath}
\end{theorem}
For any integer $p\ge 1$, we have 
\begin{equation}
\bigg(-\frac{\log_{\lambda}(1-x)}{1-x}\bigg)^{p}=\bigg(\sum_{k=1}^{\infty}H_{k,\lambda}x^{k}\bigg)^{p}=\sum_{k=p}^{\infty}\bigg(\sum_{k_{1}+\cdots+k_{p}=k}H_{k_{1},\lambda}\cdots H_{k_{p},\lambda}\bigg)x^{k}.\label{39}
\end{equation}
To ease notation, we let 
\begin{equation}
H_{k,p,\lambda}=\sum_{k_{1}+\cdots+k_{p}=k}H_{k_{1},\lambda}H_{k_{2},\lambda}\cdots H_{k_{p},\lambda}, \quad(k \ge p), \label{40}
\end{equation}
where the sum runs over all positive integers $k_{1}, \dots, k_{p}$, with $k_{1}+\cdots+k_{p}=k$.
Then, by \eqref{39} and \eqref{40}, we get 
\begin{displaymath}
\bigg(-\frac{\log_{\lambda}(1-x)}{1-x}\bigg)^{p}=\sum_{k=p}^{\infty}H_{k,p,\lambda}x^{k}.
\end{displaymath}
For any integer $p \ge 1$, we have 
\begin{align}
&\int_{0}^{1}\bigg(\sum_{k=p}^{\infty}H_{k,p,\lambda}x^{k-1}\bigg)(1-x)^{n}dx=\sum_{k=p}^{\infty}H_{k,p,\lambda}\int_{0}^{1}x^{k-1}(1-x)^{n}dx \label{41} \\
&=\sum_{k=p}^{\infty}H_{k,p,\lambda}B(k,n+1)=n!\sum_{k=p}^{\infty}\frac{H_{k,p,\lambda}}{k(k+1)\cdots(k+n)}. \nonumber
\end{align}
From \eqref{41}, we note that 
\begin{align}
&\sum_{k=p}^{\infty}\frac{n!H_{k,p,\lambda}}{k(k+1)\cdots(k+n)}=\int_{0}^{1}\bigg(\sum_{k=p}^{\infty}H_{k,p,\lambda}x^{k}\bigg)\frac{(1-x)^{n}}{x}dx \label{42} \\
&=\int_{0}^{1}\bigg(-\frac{\log_{\lambda}(1-x)}{1-x}\bigg)^{p}\frac{(1-x)^{n}}{x}dx=\int_{0}^{1}\Big(-\log_{\lambda}(1-x)\Big)^{p}\frac{(1-x)^{n-p}}{x}dx \nonumber \\
&=\int_{0}^{1/\lambda}\frac{t^{p}e_{\lambda}^{n-p}(-t)}{1-e_{\lambda}(-t)} e_{\lambda}^{1-\lambda}(-t)dt=\int_{0}^{1/\lambda}\frac{t^{p}e^{n-p+1-\lambda}(-t)}{1-e_{\lambda}(-t)}dt \nonumber\\
&=\sum_{m=0}^{\infty}\int_{0}^{1/\lambda}t^{p}(1-\lambda t)^{\frac{n-p+1-\lambda+m}{\lambda}}dt, \nonumber
\end{align}
where we made the change of variables $ x \rightarrow 1-e_{\lambda}(-t)$.
Now, by integrating by parts we see that 
\begin{align}
\int_{0}^{1/\lambda}&t^{p}(1-\lambda t)^{\frac{n-p+1-\lambda+m}{\lambda}}dt 	\label{43} \\
&=\frac{p}{(n-p+1-\lambda+m+\lambda)}\int_{0}^{1/\lambda}t^{p-1}(1-\lambda t)^{\frac{n-p+1-\lambda+m+\lambda}{\lambda}}dt. \nonumber 
\end{align}
By repeating this process, we obtain
\begin{align*}
&\int_{0}^{1/\lambda}t^{p}(1-\lambda t)^{\frac{n-p+1-\lambda+m}{\lambda}}dt \nonumber\\
&=\frac{p(p-1)}{(n-p+1-\lambda+m+\lambda)(n-p+1-\lambda+m+2\lambda)}\int_{0}^{1/\lambda}t^{p-2}(1-\lambda t)^{\frac{n-p+1-\lambda+m+2\lambda}{\lambda}}dt \nonumber \\
&=\cdots \nonumber \\
&=\frac{p!}{(n-p+1-\lambda+m+\lambda)(n-p+1-\lambda+m+2\lambda)\cdots(n-p+1-\lambda+m+p\lambda)}\nonumber\\
&\qquad \times\int_{0}^{1/\lambda}(1-\lambda t)^{\frac{n-p+1-\lambda+m+p\lambda}{\lambda}}dt \nonumber\\ 
&=\frac{p!}{(n-p+1-\lambda+m+\lambda)(n-p+1-\lambda+m+2\lambda)\cdots(n-p+1-\lambda+m+(p+1)\lambda)}\nonumber\\
&=\frac{p!}{(n-p+1+m)(n-p+1+m+\lambda)\cdots(n-p+1+m+p\lambda)}\nonumber \\
&=p!\frac{1}{\langle n-p+1+m\rangle_{p+1,\lambda}}, \nonumber
\end{align*}
where $\langle x\rangle_{0,\lambda}=1,\ \langle x\rangle=x(x+\lambda)\cdots(x+(n-1)\lambda),\ (n\ge 1)$. \\
By \eqref{42} and \eqref{43}, we get 
\begin{equation}
\sum_{k=p}^{\infty}\frac{n!H_{k,p,\lambda}}{k(k+1)\cdots(k+n)}=p!\sum_{m=0}^{\infty}\frac{1}{\langle n-p+1+m\rangle_{p+1,\lambda}}.\label{44}
\end{equation}
Therefore, by \eqref{44}, we obtain the following theorem. 
\begin{theorem}
For any integer $p \ge 1$, we have 
\begin{displaymath}
\sum_{k=p}^{\infty}\frac{H_{k,p,\lambda}}{\binom{k+n-1}{n}(k+n)}=p!\sum_{m=0}^{\infty}\frac{1}{\langle n-p+1+m\rangle_{p+1,\lambda}},
\end{displaymath}
where $H_{k,p,\lambda}=\sum_{k_{1}+\cdots+k_{p}=k}H_{k_{1},\lambda}H_{k_{2},\lambda}\cdots H_{k_{p},\lambda}$,\,\, with $k_{1},k_{2},\dots k_{p} \ge 1$.
\end{theorem}
The degenerate Hurwitz zeta function is defined as
\begin{displaymath}
	\zeta_{\lambda}^{*}(k,x)=\sum_{m=0}^{\infty}\frac{1}{\langle m+x\rangle_{k,\lambda}}, 
\end{displaymath}
where $k\in\mathbb{N}$ with $k\ge 2$ and $x>0$. \\
Note that 
\begin{displaymath}
	\lim_{\lambda\rightarrow 0}\zeta_{\lambda}^{*}(k,x)=\sum_{m=0}^{\infty}\frac{1}{(m+x)^{k}}=\zeta(k,x),
\end{displaymath}
where $\zeta(k,x)$ is the Hurwitz zeta function. \\
From Theorem 6, we have the next corollary.
\begin{corollary}
For any integer $p \ge 1$, the following identities hold true:
\begin{align*}
&\sum_{k=p}^{\infty}\frac{H_{k,p,\lambda}}{\binom{k+n-1}{n}(n+k)}=p!\zeta_{\lambda}^{*}(p+1,n-p+1), \\
&=\sum_{m=0}^{\infty}\frac{1}{\binom{n-p+1+m+(p-1)\lambda}{p}_{\lambda}(n-p+1+m+p\lambda)}, 
\end{align*}
where 
\begin{displaymath}
\binom{x}{n}_{\lambda}=\frac{(x)_{n,\lambda}}{n!}
\end{displaymath}
are the degenerate binomial coefficients (see [15,16]).
\end{corollary}

 \section{Conclusion} 
We studied the degenerate harmonic numbers $H_{k,\lambda}$ and the degenerate higher-order harmonic numbers $H_{k,\lambda}^{(p)}$, in connection with the degenerate zeta function $\zeta_{\lambda}(s)$ and the degenerate Hurwitz zeta function $\zeta_{\lambda}^{*}(k,x)$. We showed that several infinite sums involving the degenerate higher-order harmonic numbers can be expressed in terms of the degenerate zeta functions. Some of these are as follows:
\begin{align*}
&\sum_{k=1}^{\infty}\frac{n!H_{k,\lambda}^{(p)}}{k(k+1)\cdots(k+n)}=(n-1)\sum_{j=0}^{n-2}\binom{n-2}{j}(-1)^{j} \\
&\times \left\{\sum_{k=1}^{p}\frac{(-1)^{k-1}}{(j+1)^{k}}\zeta_{\lambda}(p-k+2)+\frac{(-1)^{p}}{(j+1)^{p+1}}\sum_{i=0}^{j}\frac{i!}{(\lambda+1)(\lambda+2)\cdots(\lambda+i+1)}\right\},\\
&\sum_{k=1}^{\infty}\frac{H_{k,\lambda}^{(p)}}{k(k+1)}=\zeta_{\lambda}(p+1),\\
&\sum_{k=1}^{\infty}\frac{H_{k,\lambda}^{(p)}}{k(k+2)}=\frac{1}{2}\Big(\zeta_{\lambda}(p+1)+ \zeta_{\lambda}(p)-\zeta_{\lambda}(p-1)+\cdots+(-1)^{p} \zeta_{\lambda}(2)+\frac{(-1)^{p+1}}{\lambda+1}\bigg). 
\end{align*}\par
Furthermore, we demonstrated that an infinite sum involving finite sums of products of the degenerate harmonic numbers can be represented by using the degenerate Hurwitz zeta function. Indeed, we have shown that
\begin{equation*}
\sum_{k=p}^{\infty}\frac{n!}{k(k+1)\cdots(k+n)}\times \sum_{k_{1}+\cdots+k_{p}=k}H_{k_{1},\lambda}H_{k_{2},\lambda}\cdots H_{k_{p},\lambda}=p!\zeta_{\lambda}^{*}(p+1,n-p+1).
\end{equation*} \par
As one of our research projects, we would like to continue to study degenerate versions of some special numbers and polynomials and to find their applications to physics, science and engineering as well as to mathematics.

\end{document}